\documentclass[11pt,leqno]{article} 
\usepackage{graphics}
\newtheorem{thm}{Theorem}[section]
\newtheorem{lma}{Lemma}[section]
\newtheorem{cor}{Corollary}

\newcommand{\beqa}{\begin{eqnarray}}
\newcommand{\eeqa}{\end{eqnarray}}

\newcommand{\pf}{\noindent {\bf Proof:} $\s$ }
\newcommand{\epf}{ \hfill$\diamondsuit$ \medskip}

\newcommand{\ds}{\displaystyle}
\newcommand{\beq}{\begin{equation}}
\newcommand{\eeq}{\end{equation}}
\newcommand{\lbl}{\label}
\newcommand{\s}{\; \;}

\newcommand{\ep}{\epsilon}

\newcommand{\la}{\lambda}

\newcommand{\ra}{\rightarrow}

\newcommand{\p}{\varphi}

\title{Exact multiplicity of solutions for some semilinear Dirichlet problems}

\author{
Philip Korman   \\ 
Department of Mathematical Sciences \\ 
University of Cincinnati \\ 
Cincinnati Ohio 45221-0025 \\
}

\date{}

\begin{document}

\maketitle
\begin{abstract} 
The classical result of A. Ambrosetti and G.  Prodi \cite{AP}, in the form of M.S. Berger and E. Podolak \cite{BP},
gives the exact number of solutions for the problem 
\[
\Delta u+g(u)= \mu \phi _1(x)+e(x) \s \mbox{in $D$} , \s u=0 \s \mbox{on $\partial D$} \,,
\]
depending on the real parameter $\mu$, for a class of convex $g(u)$, and $\int _D  e(x) \phi _1(x)\, dx=0$ (where $\phi _1(x)>0$ is the principal eigenfunction  of the Laplacian on $D$, and $D \subset R^n$ is a smooth domain). By considering generalized harmonics, we give a similar result for the problem
\[
\Delta u+g(u)= \mu f(x) \s \mbox{in $D$} , \s u=0 \s \mbox{on $\partial D$} \,,
\]
with $f(x)>0$. Such problems  occur, for example,  in ``fishing" applications that we discuss, and propose a new model.
\smallskip

Our approach also produces a very simple proof of the anti-maximum principle of Ph. Cl\'{e}ment and L.A. Peletier \cite{CP}.

 \end{abstract}

\begin{flushleft}
Key words:  Global solution curves, exact number of solutions, the anti-maximum principle. 
\end{flushleft}

\begin{flushleft}
AMS subject classification: 35J61, 35J25, 92D25.
\end{flushleft}

\section{Introduction}
\setcounter{equation}{0}
\setcounter{thm}{0}
\setcounter{lma}{0}

Consider the problem
\beq
\lbl{i1}
 \Delta u+g(u)= f(x) \s \mbox{in $D$} \,, \s u=0 \s \mbox{on $\partial D$} \,,
\eeq
where $D$ is a smooth domain in $R^n$, and the functions $g(u)$ and $f(x)$ are given.
Decompose $f(x)=\mu \phi _1(x)+e(x)$, where $\phi _1(x)>0$ is the principal eigenfunction  of the Laplacian on $D$ with zero boundary condition, and $\int _D  e(x) \phi _1(x)\, dx=0$.
The classical result of A. Ambrosetti and G.  Prodi \cite{AP}, in the form of M.S. Berger and E. Podolak \cite{BP}, says that if $g(u)$ is convex and asymptotically linear at $\pm \infty$, then (under an additional restriction on the slopes of $g(u)$ at $\pm \infty$) there exists a critical $\mu _0=\mu _0(e(x))$, such that the problem (\ref{i1}) has exactly two solutions for $\mu>\mu _0$, exactly one solution if $\mu=\mu _0$,
and no solutions for $\mu<\mu _0$. However, sometimes  it is desirable to have the parameter $\mu$ in front of the entire  right hand side, and to consider the problem
\beq
\lbl{i2}
 \Delta u+g(u)= \mu f(x) \s \mbox{in $D$} \,, \s u=0 \s \mbox{on $\partial D$} \,.
\eeq
Such problems occur e.g., when one considers ``fishing" applications, see S.  Oruganti et al \cite{SS}, D.G. Costa et al \cite{C}, P. Gir\~{a}o, and H. Tehrani \cite{G}, P.M. Gir\~{a}o and M. P\'{e}rez-Llanos \cite{G1}. We present an exact multiplicity result of Berger-Podolak type for the problem (\ref{i2}), provided that $f(x)>0$ on $D$. Throughout the paper, one can easily replace the Laplacian by any uniformly elliptic operator.
\smallskip

Similar result holds for  the problem
\[
 \Delta u+g(u)= \mu f(x)+e(x) \s \mbox{in $D$} \,, \s u=0 \s \mbox{on $\partial D$} \,,
\]
with $f(x)>0$ on $D$, and $\int _D  e(x) f(x)\, dx=0$, providing a generalization of the above mentioned result of M.S. Berger and E. Podolak \cite{BP}.
\medskip

Our approach involves applying the implicit function theorem for continuation of solutions in a special way. We restrict the space of solutions by keeping   the generalized first harmonic fixed, but in return allow $\mu$ to vary. Then we  compute the direction of the turn of the solution curve, similarly to P. Korman \cite{K}. We show that there is at most  one turn in case $g(u)$ is either convex or concave. 
\medskip

The well-known anti-maximum principle of Ph. Cl\'{e}ment and L.A. Peletier \cite{CP} follows easily with this approach. 
\medskip

We apply our results to a population model with fishing. We suggest a modification of the logistic model, to admit sign-changing solutions. We argue  that one needs to consider  sign-changing solutions to get complete bifurcation diagrams.

\section{The global solution curves}
\setcounter{equation}{0}
\setcounter{thm}{0}
\setcounter{lma}{0}

We assume that $D$ is a smooth domain in $R^n$, and denote by $\la _k$ the eigenvalues of the Laplacian on $D$, with zero boundary conditions, and by $\p _k(x)$ the corresponding eigenfunctions, normalized so that $\int _D \p ^2 _k(x) \, dx=1$. It is known that $\p _1(x)>0$ is simple, and $0<\la _1 <\la _2  \leq \la _3 \leq \cdots$. We denote by $H^k(D)$ the Sobolev spaces $W^{k,2}(D)$. We shall need the following generalization of Poincare's inequality.
\begin{lma}\lbl{lma:1}
Let $u(x) \in H^1_0(D) $ be such that $\int _D u(x) f(x) \, dx =0 $, for some $f(x) \in  L^2(D)$, $f(x) \not\equiv 0$. Denote $f_1=  \int _D  f(x) \p _1(x) \, dx$, and 
\[
\ds \nu =\la _1+(\la_2-\la_1) \frac{f_1^2}{||f||_{L^2}^2} \,. 
\]
Then
\beq
\lbl{0}
\int _D |\nabla u|^2 \, dx \geq \nu  \int _D u^2 \, dx \,.
\eeq
\end{lma}

\pf
By scaling of $u(x)$, we may assume that $\int _D u^2(x)  \, dx=1$. Writing $\ds u(x) =\sum _{k=1}^{\infty} u_k \p _k(x)$, $\ds f(x) =\sum _{k=1}^{\infty} f_k \p _k(x)$, we then have
\beq
\lbl{1}
\sum _{k=1}^{\infty} u^2_k=1 \,,
\eeq
\beq
\lbl{2}
\sum _{k=1}^{\infty} u_k f_k=0 \,.
\eeq
We need to show that $\int _D |\nabla u|^2 \, dx \geq \nu $. Using (\ref{1}), we estimate
\beqa
\lbl{3}
& \int _D |\nabla u|^2 \, dx=-\int _D u \Delta u \, dx=\sum _{k=1}^{\infty} \la _k u^2_k \\ \nonumber
& =\la _1+\sum _{k=2}^{\infty} \left( \la _k-\la _1 \right) u^2_k \geq \la _1+\left( \la _2 -\la _1 \right) \sum _{k=2}^{\infty} u^2_k \,. \nonumber
\eeqa
From (\ref{2})
\beq
\lbl{4}
\left( \sum _{k=2}^{\infty} u^2_k \right)^{1/2} \left( \sum _{k=2}^{\infty} f^2_k \right)^{1/2} \geq \sum _{k=2}^{\infty} u_k f_k=|-u_1f_1|=|u_1||f_1| \,.
\eeq
Set $x=\left( \sum _{k=2}^{\infty} u^2_k \right)^{1/2}$, $f=\left( \sum _{k=2}^{\infty} f^2_k \right)^{1/2}$. In view of (\ref{1}), we get from (\ref{4})
\[
xf \geq |f_1|\sqrt{1-x^2} \,,
\]
or
\[
x^2 \geq \frac{f_1^2}{f^2+f_1^2}=\frac{f_1^2}{||f||_{L^2}^2} \,,
\]
and the proof follows from (\ref{3}).
\epf

The inequality (\ref{0}) is sharp in the following sense: when $f=\p _1$,  we have $\nu=\la _2$, and one has an equal sign in (\ref{0}) at $u=\p _2$. Clearly, $\nu \leq \la _2$, and $\nu > \la _1$ if $f(x)>0$.

\begin{lma}\lbl{lma:2}
Let $(w(x), \mu) \in H^2(D) \times R$ solve the problem
\beqa\lbl{5}
& \Delta w+a(x)w=\mu f(x) \s \mbox{in $D$} \,, \s w=0 \s \mbox{on $\partial D$}   \\ \nonumber
& \int _D w(x) f(x) \, dx =0 \,, \nonumber
\eeqa
with some $f(x) \in L^2(D)$, $f(x) \not\equiv 0$. Assume that $a(x) \in C(D)$ satisfies $a(x) <\nu$ for all $x \in D$. Then $w(x) \equiv 0$, and $\mu=0$.
\end{lma}

\pf
Multiply the equation in (\ref{5}) by $w$, and integrate. By Lemma \ref{lma:1}, we have
\[
\nu \int _D w^2 \, dx \leq  \int _D |\nabla w|^2 \, dx= \int _D a(x)w^2 \, dx <\nu \int _D w^2 \, dx  \,.
\]
Hence, $w(x) \equiv 0$, and from (\ref{5}), $\mu=0$.
\epf

\begin{lma}\lbl{lma:2.1}
Consider the problem (to find $z(x)$ and $\mu^*$)
\beqa\lbl{5.1}
& \Delta z+a(x)z=\mu^* f(x)+e(x) \s \mbox{in $D$} \,, \s w=0 \s \mbox{on $\partial D$}   \\ \nonumber
& \int _D z(x) f(x) \, dx =\xi \,, \nonumber
\eeqa
where $f(x) \in L^2(D)$ satisfies $f(x) > 0$ a.e., and $a(x) \in C(D)$ satisfies $a(x) <\nu$ for all $x \in D$. Then for any $e(x)  \in L^2(D)$, and any $\xi \in R$, the problem has a solution $(z(x),\mu ^*) \in  \left(H^2(D) \cap H^1_0(D)\right) \times R$.
\end{lma}

\pf
{\em Case 1}. Assume that the operator
\[
L[z] \equiv \Delta z+a(x)z  \, : H^2(D) \cap H^1_0(D) \ra L^2(D)
\]
is invertible. We claim that
\beq
\lbl{5.2}
\int _D L^{-1}(f(x)) f(x) \, dx \ \ne 0 \,,
\eeq
where $L^{-1}$ denotes the the inverse operator of $L[z]$. Indeed, assuming otherwise, $w(x) \equiv L^{-1}(f(x))$ is not identically zero, and it satisfies (\ref{5}), with $\mu=1$, which contradicts Lemma \ref{lma:2}. Then the solution of (\ref{5.1}) is
\[
z(x)=\mu^*  L^{-1}(f(x))+ L^{-1}(e(x)) \,,
\]
and $\mu ^*$ is chosen so that $\int _D z(x) f(x) \, dx =\xi$, which we can accomplish,  in view of (\ref{5.2}).
\medskip

\noindent
{\em Case 2}. Assume that the operator $L[z]$ is not invertible. Since $a(x) <\nu \leq \la _2$, the kernel of $L[z]$ is  one-dimensional, spanned by some $\p (x)>0$. Since  $L[z]$ is a Fredholm operator of index zero, the first equation in (\ref{5.1}) is solvable if and only if its right hand side is orthogonal to $\p(x)$. We now obtain the solution $(z(x),\mu ^*)$ of (\ref{5.1}) as follows. Choose $\mu ^*$, so that $ \int _D  \left(\mu^* f(x)+e(x) \right) \p (x)\, dx =0$. Then the first equation in (\ref{5.1}) has infinitely many solutions of the form
\[
z(x)=z_0(x)+c \p(x) \,,
\]
with some $z_0(x)$. We choose the constant $c$, so that $\int _D z(x) f(x) \, dx =\xi$.
\epf

We consider next the nonlinear problem
\beq
\lbl{6}
 \Delta u+g(u)=\mu f(x) \s \mbox{in $D$} \,, \s u=0 \s \mbox{on $\partial D$} \,.
\eeq
We shall  assume that $g(u) \in C^1(R)$, and
\beq
\lbl{7}
g(u)= \left\{
\begin{array}{ll}
\gamma_1 u+b_1(u) & \mbox{if $u < 0$} \\
\gamma_2 u+b_2(u) & \mbox{if $u \geq 0$},
\end{array}
\right. 
\eeq
with real constants $\gamma _1$, $\gamma _2$, and $b_1(u)$, $b_2(u)$ bounded for all $u \in R$. Notice that we admit the case of $\gamma _2=\gamma _1$, and in particular we allow bounded $g(u)$, in case $\gamma _2=\gamma _1=0$. We shall consider strong solutions of (\ref{6}), $u(x) \in H^2(D) \cap H^1_0(D)$.
\medskip

Any function $u(x) \in L^2(D)$ can be decomposed as 
\beq
\lbl{8}
u(x)=\xi f(x)+U(x) \,, \s \mbox{with $\int _D U(x) f(x) \, dx =0$} \,, 
\eeq
for any $f(x) \in L^2(D)$.
If $f(x)>0$ a.e., we call the constant $\xi$ the {\em generalized first harmonic} of $u(x)$. 
\medskip

We shall need the following a priori estimate.
\begin{lma}\lbl{lma:3}
Assume that $f(x) \in H^2(D)$, $f(x) > 0$ a.e., and $g(u) \in C^1(R)$ satisfies the condition (\ref{7}), and $g'(u)\leq \nu_1$, for some constant $\nu_1<\nu$. Let $u(x) \in H^2(D) \cap H^1_0(D)$ be a solution of (\ref{6}), decomposed as in (\ref{8}). Then for some positive constants $c_1$ and $c_2$
\beq
\lbl{9}
|\mu |+||U||_{H^2(D)} \leq c_1|\xi|+c_2 \,.
\eeq
\end{lma}

\pf
Using the ansatz (\ref{8}) in (\ref{6}), we have
\beq
\lbl{9.1}
\s\s\s \Delta U +\xi \Delta f+g(\xi f(x)+U)=\mu f(x) \s \mbox{in $D$} \,, \s U=0 \s \mbox{on $\partial D$} \,.
\eeq
Multiplying by $U$ and integrating, we write the result  as
\beqa\lbl{10}
&  \int _D |\nabla U|^2 \, dx-\int _D \left(\xi \Delta f \right) U \, dx \\\nonumber
& -\int _D \left[g(\xi f(x)+U)-g(\xi f(x))\right]U \, dx-\int _D g(\xi f(x))U \, dx=0 \,. \nonumber
\eeqa
Using the mean value theorem, we estimate from below the third term on the left by -$\nu _1 \int _D U^2 \,dx$. If $\xi \geq 0$, then
\[
\int _D g(\xi f(x))U \, dx=\int _D \left(\gamma_2 \xi f(x)+b_2(\xi f(x)) \right) U \, dx=\int _D b_2(\xi f(x))U \, dx  \,.
\]
Using  Lemma \ref{lma:1}, we have from (\ref{10}), for any small $\ep >0$,
\[
(\nu -\nu _1)\int _D  U^2 \, dx \leq \int _D \left(\xi \Delta f \right) U \, dx+\int _D b_2(\xi f(x))U \, dx
\]
\[
\leq \ep \int _D  U^2 \, dx+c(\ep) \xi^2 \int _D \left(\Delta f \right)^2 \, dx+\ep \int _D  U^2 \, dx+c(\ep) \,,
\]
which gives us an estimate of $\int _D  U^2 \, dx$
\[
\int _D  U^2 \, dx \leq  c_1 \xi^2+c_2 \,, \s \mbox{uniformly in $\mu$}\,,
\]
with some positive constants $c_1$, $c_2$.
In case $\xi < 0$, the same estimate follows similarly. Returning to (\ref{10}), and using (\ref{7}), we have
\beq
\lbl{11}
 \int _D \left(|\nabla U|^2+U^2 \right) \, dx \leq  c_1\xi^2+c_2 \,, \s \mbox{uniformly in $\mu$}\,.
\eeq
(Here and later on, $c_1$, $c_2$ denote possibly new positive constants.)
Then
\beq
\lbl{12}
 \int _D \left(|\nabla u|^2+u^2 \right) \, dx \leq  c_1\xi^2+c_2 \,, \s \mbox{uniformly in $\mu$}\,.
\eeq

To get an estimate of $\mu$, we now multiply (\ref{6}) by $u=\xi f+U$, and integrate
\[
\xi \mu  \int _D f^2 \, dx=- \int _D |\nabla u|^2 \, dx+\int _D g(u)u \,dx  \,,
\]
which in view of (\ref{12}) implies that
\beq
\lbl{13}
|\xi| | \mu| \leq  c_1\xi^2+c_2 \,.
\eeq
(Observe that $|g(u)| \leq A|u|+B$, for some positive constants $A$, $B$, and for all $u$.) Fix some $\xi _0 >0$. Then for $|\xi| \geq \xi _0$, we conclude from (\ref{13})
\beq
\lbl{14-}
| \mu| \leq  c_1|\xi|+c_2 \,.
\eeq
In case $|\xi| \leq \xi _0$, we multiply (\ref{6}) by $\p _1$, and integrate to show that $|\mu| \leq  c_3$, for some $c_3>0$. We conclude that the bound (\ref{14-}) holds for all $\xi \in R$.
\medskip

We multiply (\ref{6}) by $\Delta u$, and integrate. Obtain
\[
\int _D \left(\Delta u \right)^2 \, dx+\int _D \Delta u \, g(u) \, dx=-\mu \int _D \nabla f \cdot \nabla u \, dx \,.
\]
Using the estimates (\ref{12}) and (\ref{14-}), we get
\[
\int _D \left(\Delta u \right)^2 \, dx \leq  c_1\xi^2+c_2 \,.
\]
Since $\Delta u=\Delta U+ \xi \Delta f$, we conclude that 
\[
\int _D \left(\Delta U \right)^2 \, dx \leq  c_1\xi^2+c_2 \,.
\]
By the elliptic estimates we obtain the desired bound on $||U||_{H^2(D)}$.
\epf

\begin{cor}\lbl{cor1}
In case $f(x)=\p _1(x)$, the second term on the left in (\ref{10}) vanishes, and we conclude that 
\[
||U||_{H^1(D)} \leq c \,, \s\s \mbox{uniformly in $\xi$ and $\mu$} \,,
\]
for some constant $c>0$.
\end{cor}

\begin{thm}\lbl{thm:1}
Assume that $f(x) \in H^2(D)$, $f(x) > 0$ a.e.,  and $g(u) \in C^1(R)$ satisfies the condition (\ref{7}), and we have  $g'(u)\leq \nu _1<\nu$ for all $u \in R$. Then for each $\xi \in (-\infty,\infty)$, there exists a unique $\mu$, for which the problem (\ref{6}) has a unique solution $u(x) \in H^2(D) \cap H^1_0(D)$, with the  generalized first harmonic equal to $\xi$.  The function $\mu =\phi (\xi)$ is smooth.
\end{thm}

\pf
We embed (\ref{6}) into a family of problems
\beq
\lbl{14}
\s\s\s \Delta u+\la _1 u+k \left(g(u)-\la _1 u \right)-\mu f(x)=0 \s \mbox{in $D$} \,, \s u=0 \s \mbox{on $\partial D$} \,,
\eeq
with $0 \leq k \leq 1$ ($k=1$ corresponds to  (\ref{6})).  When $(k=0,\mu=0)$ the problem has  solutions $u=a \p _1$, where $a$ is any constant. By choosing $a=a_0$, we can get the solution $u=a_0 \p _1$ of any generalized first harmonic $\xi ^0$. We now continue in $k$ the solutions of 
\beq 
\lbl{14.1}
\s\s\s\s\s  F(u,\mu,k) \equiv \Delta u+\la _1 u+k \left(g(u)-\la _1 u \right)-\mu f(x)=0 \; \mbox{in $D$} \,, \s u=0 \; \mbox{on $\partial D$} \eeq
\[
 \int _D uf \, dx=\xi ^0 \,, \
\]
with the operator $F(u,\mu,k) \, : H^2(D) \times R \times R \ra L^2(D)$.
We will show that the implicit function theorem applies, allowing us to continue $(u,\mu)$ as a function of $k$. Compute the Frechet derivative
\beqa \nonumber
& F_{(u,\mu)}(u,\mu, k)(w, \mu^*)=\Delta w+\la _1 w+k\left(g'(u)-\la _1 \right)w-\mu^*f(x) \,, \\ \nonumber
& \int _D wf \, dx=0 \,. \nonumber
\eeqa
By Lemma \ref{lma:2}, the map $F_{(u,\mu)}(u,\mu, k)(w, \mu^*)$ is injective, and by Lemma \ref{lma:2.1} this map is surjective. Hence, the implicit function theorem applies, and we have a solution curve $(u,\mu)(k)$. By the a priori estimate of Lemma \ref{lma:3}, this curve continues for all $0 \leq k \leq 1$, and at $k=1$, we obtain a solution of the problem (\ref{6}) with the  generalized first harmonic equal to $\xi ^0$.
\medskip

Turning to the uniqueness, let $(\bar \mu,\bar u(x))$ be another solution of (\ref{6}), and $\bar u(x)$ has the generalized first harmonic equal to $\xi ^0$. Then $(\bar \mu,\bar u(x))$ is solution of (\ref{14.1}) at $k=1$. We continue this solution backward in $k$, until $k=0$, using  the implicit function  theorem. By the Fredholm alternative, we have  $\mu =0$, when $k=0$. Then $u=a_1 \p _1$, with $a_1 \ne a_0$ (since the solution curves do not intersect), and $\bar u(x)$ has  the generalized first harmonic equal to $\xi ^0$, a contradiction.
\medskip

Finally, we show that  solutions of (\ref{6}) can be continued in $\xi$, by using the implicit function theorem. Decomposing $u(x)=\xi f(x)+U(x)$, with $\int _D Uf \, dx=0$, we see that $U(x)$ satisfies
\beqa \nonumber
& F(U,\mu,\xi) \equiv \Delta U+ g\left(\xi f(x)+U(x) \right)=\mu f(x)-\xi \Delta f \s \mbox{in $D$} \,, \s U=0 \s \mbox{on $\partial D$}
\\ \nonumber
& \int _D Uf \, dx=0 \,. \nonumber
\eeqa
Compute the Frechet derivative
\beqa \nonumber
& F_{(U,\mu)}(U,\mu, \xi)(w, \mu^*)=\Delta w+g'\left(\xi f(x)+U(x) \right)w-\mu^*f(x) \,, \\ \nonumber
& \int _D wf \, dx=0 \,. \nonumber
\eeqa
As before, we see that  the implicit function theorem applies, and we have a smooth solution curve $(u,\mu)(\xi)$ for the problem  (\ref{6}). By Lemma \ref{lma:3}, this curve continues for all $\xi \in R$.
\epf

\noindent
{\bf Remark} The theorem implies that the value of $\xi$ is a {\em global parameter}, uniquely identifying the solution pair $(\mu, u(x))$.
\medskip

The well-known anti-maximum principle is easily proved by a similar argument. As in J. Shi \cite{S}, we state it  along with the classical maximum principle. We present a self-contained proof, since the a priori estimate of Lemma \ref{lma:3} is not needed for this local result. 

\begin{thm}
Consider the following problem, with $f(x) \in L^2(D)$, and  $f(x)>0$ a.e. in $D$,
\beq
\lbl{16}
 \Delta u+\la u= f(x) \s \mbox{in $D$} \,, \s u=0 \s \mbox{on $\partial D$} \,.
\eeq
Then there exists a constant
$\delta _f$, which depends on
$f$, such that if $\la_1<\la<\la_1 +\delta _f$, then
\beq
\lbl{17}
u(x)>0  \,, \s x \in D \,, \s \frac{\partial u}{\partial n} <0 \,, \s x \in \partial D \,;
\eeq
and if $\la <\la _1$, then 
\[
u(x)<0  \,, \s x \in D \,, \s \frac{\partial u}{\partial n} >0 \,, \s x \in \partial D \,.
\]
\end{thm}

\pf
We prove the first part. Consider the problem
\beq
\lbl{18}
\Delta u+\la _1 u+k u=\mu f(x) \s \mbox{in $D$} \,, \s u=0 \s \mbox{on $\partial D$} \,.
\eeq
When $k=0$, and $\mu=0$, this problem has  a solution $u= \p _1$. Decompose $\p _1=\xi ^0 f(x)+e(x)$, with $\int _D f(x) e(x) \, dx=0$.
We now continue in $k$, $k \geq 0$ the solution $(u,\mu) \in \left(H^2(D) \cap H^1_0(D)\right) \times R$ of 
\beqa \nonumber
& \Delta u+\la _1 u+k u=\mu f(x) \s \mbox{in $D$} \,, \s u=0 \s \mbox{on $\partial D$} \\ \nonumber
& \int _D uf \, dx=\xi ^0 \,, \nonumber
\eeqa
beginning with $(\p _1,0)$ at $k=0$.
By Lemmas \ref{lma:2} and \ref{lma:2.1}, the implicit function theorem applies, and we have a solution curve $(u,\mu)(k)$, at least for small $k$. If $k>0$ is small, then $u(x)$ is close to $\p _1(x)$, and we have $u(x)>0$ in $D$, and $\frac{\partial u}{\partial n} <0$ on $\partial D$ (a.e.). Multiplying (\ref{18}) by $\p _1(x)$, and integrating over $D$, we conclude that $\mu=\mu(k)>0$. 
Then $\frac{u(x)}{\mu}$ is the solution of (\ref{16}), satisfying (\ref{17}) (a.e.).
\epf

We now study the {\em global solution curve} $\mu =\phi (\xi)$ for the problem (\ref{6}), with $u(x) \in H^2(D) \cap H^1_0(D)$, in case $g(u)$ is either convex or concave.

\begin{thm}\lbl{thm:2}
Assume that $f(x) \in H^2(D)$, $f(x) > 0$ a.e.,  and $g(u) \in C^2(R)$ satisfies the condition (\ref{7}), and we have $g'(u)\leq \nu _1<\nu$ for all $u \in R$. Assume that either $g''(u)>0$, or $g''(u)<0$ holds for all $u \in R$. Then the solution curve of  the problem (\ref{6}) $\mu =\phi (\xi)$ is either monotone, or it has exactly one critical point, which is the point of global minimum in case  $g''(u)>0$ for all $u \in R$, and the point of global maximum in case  $g''(u)<0$ for all $u \in R$.
\end{thm}

\pf
By the Theorem \ref{thm:1}, the problem (\ref{6}) has a solution curve $(u,\mu)(\xi)$, where $\xi$ is the generalized first harmonic of $u(\xi)$.
Differentiate the equation (\ref{6}) in $\xi$
\beq
\lbl{19}
 \Delta u_{\xi}+g'(u)u_{\xi}=\mu '(\xi) f(x) \s \mbox{in $D$} \,, \s u_{\xi}=0 \s \mbox{on $\partial D$} \,.
\eeq
We claim that $u_{\xi}(x) \not \equiv 0$ for all $\xi \in R$. Indeed, since $u(x)=\xi f(x)+U(x)$, we have $u_{\xi}(x)=f(x)+U_{\xi}(x)$. If $u_{\xi} (x) \equiv 0$, then $U_{\xi}(x)=-f(x)$, but $\int _D U_{\xi}(x)f(x) \,dx=0$, a contradiction.
\medskip

Assume that $\mu '(\xi_0)=0$ at some $\xi_0$. Denoting $w(x)=u_{\xi} $ at $\xi=\xi _0$, we see that $w(x)$ is a non-trivial solution of 
\beq
\lbl{20}
 \Delta w+g'(u)w=0 \s \mbox{in $D$} \,, \s w=0 \s \mbox{on $\partial D$} \,.
\eeq
Since $g'(u)<\la_2$, it follows that $w(x)>0$ in $D$. In the spirit of \cite{KLO} and \cite{OS}, we differentiate the equation (\ref{19}) once more in  $\xi$, and set $\xi=\xi _0$:
\beq
\lbl{21}
\s\s \Delta u_{\xi \xi}+g'(u)u_{\xi\xi}+g''(u)w^2=\mu ''(\xi_0) f(x) \s \mbox{in $D$} \,, \s u_{\xi \xi}=0 \s \mbox{on $\partial D$} \,.
\eeq
Combining the equations (\ref{20}) and (\ref{21}), we have
\[
\mu ''(\xi_0) \int _D w f(x) \, dx=\int _D g''(u)w^3 \, dx \,.
\]
It follows that $\mu ''(\xi_0)>0$ ($\mu ''(\xi_0)<0$) in case $g''(u)>0$ for all $u \in R$ ($g''(u)<0$ for all $u \in R$), so that any critical point of $\mu (\xi)$ is a local minimum (maximum), and hence at most one critical point is possible.
\epf

It is now easy to classify all of the possibilities.

\begin{thm}\lbl{thm:3}
Assume that $f(x) \in H^2(D)$, $f(x) > 0$ a.e.,  and $g(u) \in C^2(R)$ satisfies the condition (\ref{7}), and we have $g'(u)\leq \nu _1<\nu$ for all $u \in R$. Assume also that $g''(u)>0$ for all $u \in R$. \newline
(i) If $\gamma _1 \,, \gamma _2<\la _1$, then the problem (\ref{6}) has a unique solution for any $\mu \in R$. Moreover, the solution curve $\mu =\phi (\xi)$ is defined, and monotone decreasing for all $\xi \in R$.\newline
(ii) If $\la _1 < \gamma _1 \,, \gamma _2<\nu$, then the problem (\ref{6}) has a unique solution for any $\mu \in R$. Moreover, the solution curve $\mu =\phi (\xi)$ is defined, and  monotone increasing for all $\xi \in R$.\newline
(iii) If $  \gamma _1<\la _1 < \gamma _2<\nu$, then there is a critical $\mu _0$, so that the problem (\ref{6}) has exactly two solutions for $\mu >\mu _0$, it has a unique solution at $\mu =\mu _0$, and no solutions for $\mu <\mu _0$. Moreover, the solution curve $\mu =\phi (\xi)$  is defined  for all $\xi \in R$, it is parabola-like, and $\mu _0$ is its global minimum value.
\end{thm}

\pf
The convexity of $g(u)$ implies that $\gamma _1 <\gamma _2$. By the Theorem \ref{thm:2}, the problem (\ref{6}) has a  solution curve $\mu =\phi (\xi)$, defined for all $\xi \in R$, which  is either monotone, or it has exactly one critical point, which is the point of global minimum. 
Decompose $u(x)=\bar \xi \p_1(x)+\bar U(x)$, where $\int _D \bar U(x) \p_1(x) \, dx=0$, and $\bar \xi$ is the first harmonic. We have
\[
u(x)=\xi f(x)+ U(x)=\bar \xi \p_1(x)+\bar U(x) \,.
\]
Multiplying this by $f(x)$, and integrating
\[
\xi \int _D f^2(x) \, dx=\bar \xi \int _D f(x) \p_1(x) \, dx+\int _D \bar U(x) f(x) \, dx \,.
\]
By the Corollary \ref{cor1} of Lemma \ref{lma:3}, $\int _D \bar U(x) f(x) \, dx$ is uniformly bounded. 
It follows that $\bar \xi \ra \infty$ ($-\infty$) if an only if $\xi \ra \infty$ ($-\infty$), providing us with a ``bridge" between $\xi$ and $\bar \xi$.
\medskip
 
Multiply the equation in  (\ref{6}) by $\phi _1$, and integrate:
\beq
\lbl{22}
\mu \int _D f(x) \phi _1(x) \, dx=-\la _1 \bar \xi+\int _D   g(u)  \phi _1(x) \, dx \,,
\eeq
with $\int _D f(x) \phi _1(x) \, dx>0$. 
If $\bar \xi >0$ and large, then $u(x)=\bar \xi \p_1(x)+\bar U(x)>0$ a.e. in $D$, and  we have $u(x)<0$ a.e. in $D$ if  $\bar \xi <0$ and $|\bar \xi|$ is large.
By  the condition (\ref{7}),
\[
\mu \sim \left( \gamma _2-\la _1 \right)\bar \xi \,, \s \mbox{when $\bar \xi >0$ and large} \,,
\]
\[
\mu \sim \left( \gamma _1-\la _1 \right)\bar \xi \,, \s \mbox{when $\bar \xi<0$ and $|\bar \xi|$ is large} \,.
\]
These formulas give us the behavior of $\mu =\phi (\xi)$, as $\xi \ra \pm \infty$, and the theorem follows.
\epf

The following result is proved similarly (the concavity of $g(u)$ implies that $\gamma _2 <\gamma _1$).
\begin{thm}\lbl{thm:4}
Assume that $f(x) \in H^2(D)$, $f(x) > 0$ a.e.,  and $g(u) \in C^2(R)$ satisfies the condition (\ref{7}), and we have  $g'(u)\leq \nu _1<\nu$ for all $u \in R$. Assume also that $g''(u)<0$ for all $u \in R$. \newline
(i) If $\gamma _1 \,, \gamma _2<\la _1$, then the problem (\ref{6}) has a unique solution for any $\mu \in R$. Moreover, the solution curve $\mu =\phi (\xi)$ is monotone decreasing for all $\xi \in R$.\newline
(ii) If $\la _1 < \gamma _1 \,, \gamma _2<\nu$, then the problem (\ref{6}) has a unique solution for any $\mu \in R$. Moreover, the solution curve $\mu =\phi (\xi)$ is monotone increasing for all $\xi \in R$.\newline
(iii) If $  \gamma _2<\la _1 < \gamma _1<\nu$, then there is a critical $\mu _0$, so that the problem (\ref{6}) has exactly two solutions for $\mu <\mu _0$, it has a unique solution at $\mu =\mu _0$, and no solutions for $\mu >\mu _0$. Moreover, the solution curve $\mu =\phi (\xi)$ is parabola-like, and $\mu _0$ is its global maximum value.
\end{thm}

It appears that there is less interest in concave nonlinearities, compared with the convex ones. This may be due to the fact that if one considers {\em positive } solutions of 
\[
\Delta u+g(u)=0\s \mbox{in $D$} \,, \s u=0 \s \mbox{on $\partial D$} \,,
\]
and $g(0) \geq 0$, then the case of concave $g(u)$ is easy, and the convex case is interesting. However, if $g(0)<0$, the situation may be  reversed even for positive  solutions, see e.g., \cite{K1}. For sign-changing solutions, it seems   that the convex and concave cases are of equal complexity. 
\medskip

Examining the proofs, we see that the Theorems \ref{thm:3} and \ref{thm:4} hold verbatim for the problem
\[
\Delta u +g(u)=\mu f(x) +e(x) \s \mbox{in $D$} \,, \s u=0 \s \mbox{on $\partial D$} \,,
\]
with $e(x) \in H^2(D)$ satisfying $\int _D e(x) f(x) \, dx=0$, giving a generalization of the classical results of A. Ambrosetti and G.  Prodi \cite{AP},   and  of M.S. Berger and E. Podolak \cite{BP}.
\section{A population model with fishing}
\setcounter{equation}{0}
\setcounter{thm}{0}
\setcounter{lma}{0}

One usually begins population modeling with a logistic model 
\beq
\lbl{f1}
u'(t)=au(t)-bu^2(t) \,.
\eeq
Here $u(t)$ can be thought of  as the number of fish in a lake at time $t$; $a$ and $b$ are positive constants. When $u(t)$ is small, $u^2(t)$ is negligible, and the population grows exponentially, but after some time the growth rate decreases. Now suppose the lake occupies some region $D \subset R^n$, and $u=u(x,t)$, with $x \in D$. Suppose that fish diffuses around the lake, and the population is near zero at the banks. Assume also there is time-independent fishing, accounted by the term $\mu f(x)$, where $f(x)$ is a positive function, and $\mu $ is a parameter. Then the model is 
\[
 u_t(x,t)=au(x,t)-bu^2(x,t)+\Delta u(x,t)-\mu f(x) \s \mbox{in $D$}, \s u=0 \s \mbox{on $\partial D$} \,.
\]
We shall consider its steady state $u=u(x)$, satisfying
\beq
\lbl{f2}
\Delta u(x)+au(x)-bu^2(x)-\mu f(x) =0 \s \mbox{in $D$}, \s u=0 \s \mbox{on $\partial D$} \,.
\eeq

It is customary in the population modeling to look for positive solutions. However one does not expect the solutions of (\ref{f2}) to remain positive, when the   parameter $\mu >0 $ is varied (since $f(x)>0$). Therefore, we shall admit sign-changing solutions, with the interpretation  that some re-stocking of fish is necessary when $u(x)<0$ (which presumably occurs near the banks, i.e., $\partial D$), to avoid the algae growth or other negative consequences. However, there is no reason to use the logistic model (\ref{f1}) for sign-changing $u$. When $u<0$, it is still reasonable to assume that $u'(t) \approx a u(t)<0$, which corresponds to the assumption that the situation further deteriorates without re-stocking, but there seems to be no justification for the $-bu^2$ term.
\medskip

We  consider the following model ($f(x)>0$)
\beq
\lbl{f3}
\Delta u(x)+g(u(x))-\mu f(x) =0 \s \mbox{in $D$}, \s u=0 \s \mbox{on $\partial D$} \,,
\eeq
where $g(u)$ is an extension of the logistic model to $u<0$, which we describe next. Namely, we assume that $g(u) \in C^2(R)$, and it satisfies
\beq
\lbl{f4}
g(u)=au-bu^2 \s \mbox{for $u \geq 0$} \,, \s \mbox{with $\la _1<a<\nu$, and $b>0$} \,,
\eeq
\beq
\lbl{f5}
g'(u) <\nu \,, \s \mbox{and} \s g''(u)<0 \s \mbox{for $u \in R$} \,,
\eeq
where $ \la _1<\nu \leq \la _2$ was defined in Lemma \ref{lma:1}. Our conditions imply that $g(u) \sim cu+d$ as $u \ra -\infty$, for some constants $0<c<\nu$, and $d>0$.
\medskip

When $\mu =0$ (no fishing), the problem (\ref{f3}) has the trivial solution $u(x) \equiv 0$, and a unique positive solution $u_0(x)$, see e.g., P.  Korman and A. Leung \cite{KL}. When $\mu >0$ is varied these two solutions turn out to be connected by a smooth solution curve. To prove this result, we shall need the following consequence of Lemma 3.3 in \cite{A}.

\vspace{1.3in}

\begin{picture}(0,40)(-70,0)
\scalebox{0.85}{

\put(60, 0){\vector(1,0){130}}
\put(60,0){\vector(0,1){130}}
\multiput(106,0)(0,7){8}{\line(0,1){5}}

\thicklines

\put(44,101){\makebox(0,0)[l]{{\bf $u_0$}}}
\put(196,-1){\makebox(0,0)[l]{${\bf {\Large \mu}}$}}
\put(65,135){\makebox(0,0)[l]{${\bf ||u||}$}}
\qbezier(60,0)(150,60)(60,102)
\put(104,-8){\makebox(0,0)[l]{${\bf { \bar \mu}}$}}
%\qbezier(60,139)(165,100)(220,102)
%\qbezier(60,23)(100,26)(150,70)
%\qbezier(150,70)(190,98)(225,98)

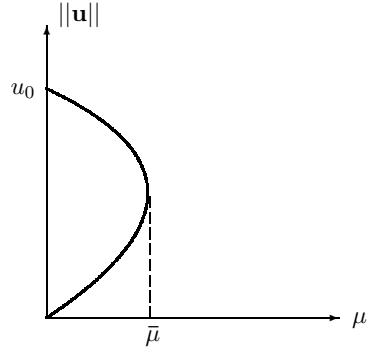
\put(10,-46){\makebox(0,0)[l]{Figure 1. The solution curve for the fishing model}}
}
\end{picture}
\vspace{.9in}

\begin{lma}\lbl{lma:A}
Let $u(x,\mu)$ denote the classical solution of (\ref{f3}), depending on a parameter $\mu \in R$. Assume that $\mu _2>\mu _1$, and  $u(x,\mu_1)>0$, $u(x,\mu_2)>0$ for all $x \in D$. Then  $u(x,\mu_2)>u(x,\mu_1)$ for all $x \in D$.
\end{lma}

Let $\xi _0=\int_D u_0(x) f(x) \, dx>0$ denote the first generalized harmonic of  $u_0(x)$. We have the following result, for possibly sign-changing solutions.

\begin{thm}\lbl{thm:f}
Assume that the conditions (\ref{f4}) and (\ref{f5}) hold, and $f(x) \in C^{\alpha}(D)$, $f(x)>0$ in $D$, $\alpha>0$.
Then in the $(\xi, \mu)$ plane there is a smooth parabola-like solution curve $\mu=\p (\xi)$ of (\ref{f3}), connecting the points $(0,0)$ and $(\xi _0,0)$. It has a unique point of maximum at some $\bar \xi \in (0,\xi _0)$, with $\bar \mu=\p (\bar \xi)>0$. Equivalently, for $\mu \in [0,\bar \mu)$ the problem  (\ref{f3}) has exactly two solutions, it has exactly one solution at $\mu =\bar \mu$, and no solutions for $\mu >\bar \mu$. Moreover, all solutions lie on a  parabola-like solution curve in the $(\mu, ||u||)$ plane, with a turn to the left.
\end{thm}

\pf
By Theorem \ref{thm:1}, we continue the solution curve from the  point $(\xi _0,0)$ in the $(\xi, \mu)$ plane for decreasing $\xi$. By Lemma \ref{lma:A}, it follows that $\mu >0$ for $\xi$ near $\xi _0$, and  $\xi<\xi _0$ (if $\mu <0$, then $\xi >\xi _0$).  By Theorem \ref{thm:4}, this curve $\mu=\mu(\xi)$ has a unique critical point on $ (0,\xi _0)$, which a point of global maximum, and this curve links up to the point $(0,0)$. This implies that the solution curve  in the $(\mu, ||u||)$ plane is as in Figure 1, concluding the proof.
\epf

S.  Oruganti et al \cite{SS} considered positive solutions of (\ref{f2}). They proved a similar result (as in Figure 1) for $a$ sufficiently close to $\la _1$. In that case, $\xi _0$ is small, and the entire solution curve is close to the point $(0,0)$. Working with positive  solutions only  narrows the class of solutions considerably, and the result of \cite{SS} is probably the best one can get (for the picture as in Figure 1). We showed in \cite{K2} that the picture is  different when $a>\la _2$ (in case of positive solutions).

\vspace{1.3in}

\begin{picture}(0,40)(-70,0)
\scalebox{0.85}{

\put(20, 0){\vector(1,0){170}}
\put(60,0){\vector(0,1){130}}
\multiput(106,0)(0,7){7}{\line(0,1){5}}

\thicklines

\put(66,103){\makebox(0,0)[l]{{\bf $u_0$}}}
\put(196,-1){\makebox(0,0)[l]{${\bf {\Large \mu}}$}}
\put(65,135){\makebox(0,0)[l]{${\bf ||u||}$}}
\qbezier(60,0)(170,60)(10,125)
\put(104,-8){\makebox(0,0)[l]{${\bf { \bar \mu}}$}}
%\qbezier(60,139)(165,100)(220,102)
%\qbezier(60,23)(100,26)(150,70)
%\qbezier(150,70)(190,98)(225,98)

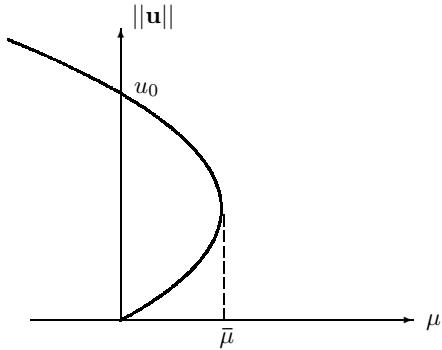
\put(-20,-46){\makebox(0,0)[l]{Figure 2. The  fishing model with stocking of fish (when $\mu <0$)}}
}
\end{picture}
\vspace{.9in}

We show next that the upper branch of the solution curve in Theorem \ref{thm:f}  continues for $\xi \in (-\infty,0)$, with     $\mu=\p (\xi)$ monotone decreasing, and $u(x)>0$ in $D$. Moreover, $\lim _{\xi \ra -\infty} \p (\xi))=+\infty$,  implying that   the solution curve in the $(\mu, ||u||)$ plane is as in  Figure 2. Indeed, let us return to the solution point $(\xi _0,0)$ in the $(\xi, \mu)$ plane. For $\xi >\xi _0$, we have $\mu <0$ by Lemma \ref{lma:A}. Then $u(x)>0$, by the minimum principle, so that $g'(u)=a-2u<\nu$, and Theorem \ref{thm:4} applies. One can interpret $\mu<0$ as {\em stocking of fish}.

\end{document}